\documentclass[10pt]{article}

\usepackage{amsmath,amscd,amsfonts,amsthm}

\newcommand{\shrinkmargins}[1]{
  \addtolength{\textheight}{#1\topmargin}
  \addtolength{\textheight}{#1\topmargin}
  \addtolength{\textwidth}{#1\oddsidemargin}
  \addtolength{\textwidth}{#1\evensidemargin}
  \addtolength{\topmargin}{-#1\topmargin}
  \addtolength{\oddsidemargin}{-#1\oddsidemargin}
  \addtolength{\evensidemargin}{-#1\evensidemargin}
  }

\shrinkmargins{.7}

\DeclareMathOperator{\SL}{SL}

\DeclareMathOperator{\PSL}{PSL}
\DeclareMathOperator{\GL}{GL}

\DeclareMathOperator{\End}{End}

\DeclareMathOperator{\Spec}{Spec}

\DeclareMathOperator{\Gal}{Gal}

\newcommand{\field}[1]{\mathbb{#1}}

\newcommand{\Qp}{\field{Q}_p}
\newcommand{\Q}{\field{Q}}

\newcommand{\Z}{\field{Z}}
\newcommand{\F}{\field{F}}

\newcommand{\fivebar}{\bar{\field{F}}_5}
\newcommand{\R}{\field{R}}
\newcommand{\C}{\field{C}}

\renewcommand{\P}{\field{P}}
\newcommand{\ord}{\mbox{ord}}

\renewcommand{\AA} {\mathcal{A}}
\newcommand{\XX} {\mathcal{X}}
\newcommand{\ra}{\rightarrow}

\newcommand{\OO}{\mathcal{O}}

\newcommand{\ic}[1]{\mathfrak{#1}}

\newcommand{\mat}[4]{\left[\begin{array}{cc}#1 & #2 \\
                                         #3 & #4\end{array}\right]}

\newcommand{\inj}{\hookrightarrow}
\newcommand{\tensor} {\otimes}

\newcommand{\GalK}{\Gal(\bar{K}/K)}
\newcommand{\GalF}{\Gal(\bar{F}/F)}
\newcommand{\GalQ}{\Gal(\bar{\Q}/\Q)}

\newcommand{\beq}{\begin{displaymath}}
\newcommand{\eeq}{\end{displaymath}}
\newcommand{\beqn}{\begin{equation}}
\newcommand{\eeqn}{\end{equation}}

\theoremstyle{plain}
\newtheorem{thm}{Theorem}[section]
\newtheorem{prop}[thm]{Proposition}
\newtheorem{cor}[thm]{Corollary}
\newtheorem{lem}[thm]{Lemma}
\newtheorem*{intro}{Theorem}
\newtheorem*{introconj}{Conjecture}
\newtheorem*{introcor}{Corollary}

\theoremstyle{definition}
\newtheorem{defn}[thm]{Definition}

\theoremstyle{remark}

\title{Serre's conjecture over $\F_9$}

\author{Jordan S. Ellenberg}
\date{10 May 2001}

\begin{document}

\maketitle

\section*{Abstract}

In this paper we show that an odd Galois representation $\bar{\rho}:\GalQ
\ra \GL_2(\F_9)$ satisfying certain local conditions at $3$ and $5$ is
modular.  Our main tool is an idea of Taylor~\cite{tayl:icos2}, which
reduces the problem to that of exhibiting points on a Hilbert modular
surface which are defined over a solvable extension of $\Q$, and which
satisfy certain reduction properties.  As a corollary, we show that
Hilbert-Blumenthal abelian surfaces with good ordinary reduction at
$3$ and $5$ are modular.

\section*{Introduction}

In 1986, J.~P.~Serre proposed the following conjecture:

\begin{introconj} Let $\F$ be a finite field of characteristic $p$, and 
\beq
\bar{\rho}: \GalQ \ra \GL_2(\F)
\eeq
an irreducible representation such that $\det \bar{\rho}$ applied to
complex conjugation yields $-1$.  Then $\bar{\rho}$ is the mod $p$
representation attached to a modular form on $GL_2(\Q)$.
\end{introconj}

Serre's conjecture, if true, would provide the first serious glimpse
into the non-abelian structure of $\GalQ$.  The work of Langlands and
Tunnell shows that Serre's conjecture is true when $\GL_2(\F)$ is
solvable; that is, when $\F$ is $\F_2$ or $\F_3$.  Work of
Shepherd-Barron and Taylor~\cite{shep:icos} and
Taylor~\cite{tayl:icos2} have shown that the conjecture is also true,
under some local conditions on $\bar{\rho}$, when $\F$ is $\F_4$ or
$\F_5$.

In the present work, we show that Serre's conjecture is true, again
subject to certain local conditions, when $\F = \F_9$.  To be precise,
we prove the following theorem.

\begin{intro}
Let
\beq
\bar{\rho}: \GalQ \ra \GL_2(\F_9)
\eeq
be an odd Galois representation such that
\begin{itemize}
\item  The restriction of $\bar{\rho}$ to $D_3$ can be written as
\beq
\bar{\rho}|D_3 \cong \mat{\psi_1}{*}{0}{\psi_2},
\eeq
where $\psi_1$ and $\psi_2$ are characters such that $(\psi_1 \oplus
\psi_2)|I_3$ is isomorphic to the direct sum of the trivial character
and $\bar{\chi}_3$;

\item The image of the inertia group $I_5$ lies in $\SL_2(\F_9)$, and
has odd order in $\PSL_2(\F_9)$.

\end{itemize}

Then $\bar{\rho}$ is modular.
\end{intro}

As a corollary, we get the following result towards a generalized
Shimura-Taniyama-Weil conjecture for abelian surfaces:

\begin{introcor}  Let $A/\Q$ be a Hilbert-Blumenthal abelian surface which has
good ordinary reduction at $3$ and $5$.  Then $A$ is a quotient of
$J_0(N)$ for some integer $N$. 
\end{introcor}

We prove the theorem above by exhibiting $\bar{\rho}$ as the Galois
representation on the $3$-torsion subscheme of a certain
Hilbert-Blumenthal abelian surface defined over a totally real
extension $F/\Q$ with solvable Galois group.  We then use an idea of
Taylor, together with a new theorem of Skinner and
Wiles~\cite{skin:sw2}, to prove the modularity of the abelian surface,
and consequently of $\bar{\rho}$.

The key algebro-geometric point is that a certain twisted Hilbert
modular variety has many points defined over solvable extensions of
$\Q$.  This suggests that we consider the class of varities $X$ such
that, if $K$ is a number field, and $\Sigma$ is the set of all
solvable Galois extensions $L/K$, then
\beq
\bigcup_{L \in \Sigma} X(L)
\eeq
is Zariski-dense in $X$.  We say $X$ has ``property {\bf S}'' in this
case.  Certainly if $X$ has a Zariski-dense set of points over a {\em
single} number field---for example, if $X$ is unirational---it has
property {\bf S}.  The Hilbert modular surfaces we consider, on the
other hand, are varieties of general type with property {\bf S}.

To indicate our lack of knowledge about solvable points on varieties,
note that at present there does not exist a variety which we can prove
does {\em not} have property {\bf S}!  Nonetheless, it seems reasonable
to guess that ``sufficiently complicated'' varieties do not have
property {\bf S}. 

\medskip

One might consider the present result evidence for the
truth of Serre's conjecture.  On the other hand, it should be pointed out
that the theorems here and in \cite{shep:icos}, \cite{tayl:icos2} rely
crucially on the facts that
\begin{itemize}
\item the $\GL_2$ of small finite fields is solvable, and
\item certain Hilbert modular varieties for number fields of small
discriminant have property {\bf S}.
\end{itemize}

These happy circumstances may not persist very far.  In particular, it
is reasonable to guess that only finitely many Hilbert modular
varieties have property {\bf S}.  If so, one might say that we have
much philosophical but little numerical evidence for the truth of
Serre's conjecture in general.  Our ability to compute
has progressed mightily since Serre's conjecture was first announced.
It would be interesting, given the present status of the conjecture,
to carry out numerical experiments for $\F$ a ``reasonably large''
finite field---whatever that might mean.

\section*{Notation}

If $\bar{\rho}: \GalK \ra \GL_2(\F_9)$ is a Galois representation, we
define $V_{\bar{\rho}}$ to be the symplectic Galois module $\F_9
\oplus \F_9$, with Galois acting by $\bar{\rho}$, and the standard
symplectic pairing.

If $K$ is a totally real number field, 
write $c \in \GalK$ for a complex conjugation.  If $v$ is a prime of a
number field $F$, we write $D_v \subset G_F$ for the decomposition
group associated to $v$, and $I_v$ for the corresponding inertia
group.  The $p$-adic cyclotomic character of Galois is denoted by
$\chi_p$, and its mod $p$ reduction by $\bar{\chi}_p$.

If $V \subset \P^N$ is a projective variety, write $F_1(V)$ for the
Fano variety of lines contained in $V$.

If $\OO$ is a ring, an {\em $\OO$-module scheme} is an $\OO$-module in
the category of schemes.

All Hilbert modular forms are understood to be of parallel weight.
 
\section{Realizations of Galois representations on HBAV's}

Recall that a {\em Hilbert-Blumenthal abelian variety} (HBAV) over a number
field is an abelian $d$-fold endowed with an injection $\OO \inj
\End(A)$, where $\OO$ is the ring of integers of a totally real number
field of degree $d$ over $\Q$.  Many Hilbert-Blumenthal abelian
varieties can be shown to be modular; for example, see~\cite{shep:icos}.
It is therefore sometimes possible to show that a certain mod $p$
Galois representation $\bar{\rho}$ is modular by realizing it on the
$p$-torsion subscheme of some HBAV.

We will show that, given a Galois representation $\GalK \ra \GL_2(\F_9)$
satisfying some local conditions at $3,5$ and $\infty$, we can find an
abelian surface over a solvable extension of $K$, satisfying some
local conditions at $3$ and $5$.  One of these conditions---that
certain representations be ``$D_p$-distinguished''---requires further comment.

\begin{defn}  Let $\bar{\rho}: \GalK \ra \GL_2(\bar{\F}_p)$ be a Galois
representation, and let $\ic{p}|p$ be a prime of $K$.  We say that
$\bar{\rho}$ is {\em $D_\ic{p}$-distinguished} if the
semisimplification of the restriction $\GalK|D_\ic{p}$ is isomorphic
to $\theta_1 \oplus \theta_2$, with $\theta_1$ and $\theta_2$ distinct
characters from $\GalK$ to $\bar{\F}^*_p$.
\end{defn}

This condition is useful in deformation theory, and is required, in
particular, in the main theorem of \cite{skin:sw2}.  A natural source
of $D_\ic{p}$-distinguished Galois representations is provided by
abelian varieties with ordinary reduction at $\ic{p}$.

\begin{prop} Let $p$ be an odd prime.  Let $K_v$ be a finite extension
of $\Q_p$ with odd ramification degree, and let $A/K_v$ be a
HBAV with good ordinary or multiplicative reduction and real
multiplication by $\OO$, and let $\ic{p}$ be a prime of $\OO$ dividing $p$.

Then the
semisimplification of the $\Gal(\bar{K}_v/K_v)$-module $A[\ic{p}]$ is
isomorphic to $\theta_1 \oplus \theta_2$, with $\theta_1$ and
$\theta_2$ distinct characters of $\Gal(\bar{K}_v/K_v)$.
\label{pr:distinguished}
\end{prop}

\begin{proof}
We think of $A[\ic{p}]$ as a two-dimensional $\OO/\ic{p}\OO$-module,
where $\OO$ is the ring of real multiplications.  Let $e$ be the
ramification index of $K_v$ over $\Q_p$.

The decomposability of $A[\ic{p}]$ into an extension of $\theta_1$ by
$\theta_2$ is a general fact about ordinary abelian varieties.

Now the action of the inertia group $I_v$ on $A[\ic{p}]$ is an
extension of the trivial character by the cyclotomic character.  Since
the ramification index of $K_v/\Q_p$ is odd, these two characters are
distinct.  Therefore, $\theta_1$ and $\theta_2$ are also distinct. 
\end{proof}

We are now ready to state the main theorem of this section.

\begin{prop}
Let $K$ be a totally real number field, and let 
\beq
\bar{\rho}: \GalK \ra \GL_2(\F_9)
\eeq
be a Galois representation such that $\det \bar{\rho} =
\bar{\chi}_3$.  Suppose that
\begin{itemize}

\item The absolute ramification degree of $K$ is odd at every prime of
$K$ above $3$ and $5$.

\item  For any prime $w$ of $K$ over $3$, the restriction of
$\bar{\rho}$ to the decomposition group $D_w$ can be written as
\beq
\bar{\rho}|D_w \cong \mat{\psi_1}{*}{0}{\psi_2},
\eeq
where $\psi_1 \oplus \psi_2 | I_w$ is isomorphic to the sum of the
trivial character and $\bar{\chi}_3$.

\item The image of the inertia group $I_v$ in $\GL_2(\F_9)$ has odd
order for every prime $v$ of $K$ over $5$.

\end{itemize}

\medskip

Then there
exists a totally real number field $F$ with $F/K$ a solvable Galois
extension, and a Hilbert-Blumenthal abelian variety $A/F$ with real
multiplication by $\OO = \OO_{\Q[\sqrt{5}]}$, such that
\begin{itemize}

\item The absolute ramification degree of $F$ is odd at every prime of
$F$ over $3$ and $5$;

\item $A$ has good ordinary or multiplicative reduction at all
primes of $F(\zeta_3)$ above $3$ and $5$;

\item $A[\sqrt{5}]$ is an absolutely irreducible $\GalF$-module;

\item there exists a symplectic isomorphism of $\GalF$-modules
\beq
\iota: A[3] \cong V_{\bar{\rho}}.
\eeq
\end{itemize}
\label{pr:hbav}
\end{prop}

\section{Proof of Proposition~\ref{pr:hbav}}

Our main tool is an expicit description of the moduli space
of HBAV's with real multiplication by $\OO$ and full $3$-level
structure, worked out by Hirzebruch and van der Geer.

\begin{lem}
Let $S_\Z$ be the surface in $\P^5/\Z$ defined by the equations
\beq
\sigma_1(s_0, \ldots, s_5)  = \sigma_2(s_0, \ldots, s_5) =
\sigma_4(s_0, \ldots, s_5) = 0,
\eeq
where $\sigma_i$ is the $i$th symmetric polynomial.  Note that $A_6 \cong
\PSL_2(\F_9)$ acts on $S_\C$ by permutation of coordinates.

Let $X_\C$ be the Hilbert modular surface parametrizing complex HBAV's with real
multiplication by $\OO$ and full $3$-level structure.  Then $X_\C$ and
$S_\C$ are related by a birational isomorphism which is compatible
with the action of $\PSL_2(\F_9)$ on $X_\C$ and $A_6$ on $S_\C$.
\label{le:vdg}
\end{lem}

\begin{proof}
\cite[VIII.(2.6)]{vdge:hms}
\end{proof}

More precisely, the complement of the cusps in $X_\C$ maps to an open
subvariety of $S_\C$.  Write $Z_\C \subset S_\C$ for the complement of the
image of this map.  We write
$Q_1, Q_2, Q_4$ for the hypersurfaces in $\P^5/\Z$ cut out by
$\sigma_1 = 0, \sigma_2 = 0, \sigma_4 =0$ respectively.

The morphism from $X_\C$ to $S_\C$ is defined as follows.  Let
$\Gamma$ be the kernel of the reduction map $\SL_2(\OO) \ra
\SL_2(\F_9)$.  Then the complex vector space $S_2(\Gamma;
\C)$  of weight
$2$ cuspforms for $\Gamma$ has dimension $5$, and $\PSL_2(\F_9)$ acts on
$S_2(\Gamma;\C)$ through its irreducible $5$-dimensional
representation~\cite[VIII,2.5]{vdge:hms}.  
Let $s_0, \ldots, s_5$ be a basis for $S_2(\Gamma;\C)$ such that
$\PSL_2(\F_9)$ acts by even permutations on $s_0, \ldots, s_5$.  Then
taking $s_0, \ldots, s_5$ as projective coordinates on $X_\C$ yields the
rational morphism from $X_\C$ to  $\P^5$ described above.

In particular, we have
\beq
\sigma_1(s_0, \ldots, s_5)  = \sigma_2(s_0, \ldots, s_5) =
\sigma_4(s_0, \ldots, s_5) = 0.
\eeq

Note that the other symmetric functions $\sigma_3, \sigma_5$, and
$\sigma_6$ in $s_0, \ldots, s_5$ are cuspforms of weight
$6,10,12$ for the whole group $\SL_2(\OO)$.  In section \ref{ss:char5}
we will described these forms in terms of classical generators of the
ring of level $1$ modular forms.  

\medskip

The complex moduli space above descends to one defined over a ring of
algebraic numbers.  Let $N$ be the product of all primes where
$\bar{\rho}$ is ramified.  (In particular, $3$ divides $N$.)  Then the
representation $\bar{\rho}$ defines an etale group scheme
$V_{\bar{\rho}}/\Z[1/N]$.  The determinant condition on $\bar{\rho}$
yields a natural isomorphism $\wedge^2 V_{\bar{\rho}} \cong \mu_3
\tensor_\Z \OO$.

  Then there exists a proper scheme
$\XX^{\bar{\rho}}/\Z[1/N]$ parametrizing pairs $(A,\iota)$, where $A$
is a generalized HBAV, and
\beq
\iota: A[3] \cong V_{\bar{\rho}}
\eeq
is an isomorphism of group schemes such that  $\wedge^2 \iota$ is the
identity map.  Write $X^{\bar{\rho}}/K$ for the restriction of
$\XX^{\bar{\rho}}$ to $\Spec K$.

The surface $X^{\bar{\rho}}/K$ is geometrically isomorphic to $X$; in
particular, its geometric genus is still $5$.  Let $t_0, \ldots, t_5$
be a set of sections spanning $H^0(X^{\bar{\rho}}, \Omega^1_{X^{\bar{\rho}}})$.  Then $t_0, \ldots, t_5$ can
be expressed as complex linear combinations of $s_0, \ldots, s_5$.  In
particular, the rational map $X^{\bar{\rho}}/K \ra \P^5/K$ defined by
$t_0, \ldots, t_4,t_5$ is
an locally closed immersion of the complement of the cusps in $X^{\bar{\rho}}$
into $\P^5$ as an open subset in an intersection of smooth
hypersurfaces of degrees $1$, $2$ and $4$. Write $S^{\bar{\rho}}/K$ for the
closure of the image of $X^{\bar{\rho}}$ in $\P^5$.

Note that there exists an isomorphism $\alpha:
\P^5_{\bar{\Q}} \ra \P^5_{\bar{\Q}}$, depending on our choice of $t_0,
\ldots, t_5$, such that $S^{\bar{\rho}} =
\alpha^{-1} S$.  Write $Q_1^{\bar{\rho}}, Q_2^{\bar{\rho}}$ and
$Q_4^{\bar{\rho}}$ for the hypersurfaces $\alpha^{-1} Q_1, \alpha^{-1} Q_2,
\alpha^{-1} Q_4$.  Then $S^{\bar{\rho}}$ is the complete intersection of
$Q_1^{\bar{\rho}}, Q_2^{\bar{\rho}}$ and $Q_4^{\bar{\rho}}$, which are
hypersurfaces in $\P^5/\Q$.

To prove Proposition~\ref{pr:hbav}, we will need to find a point on a
twisted Hilbert modular variety $X^{\bar{\rho}}$ defined over a
solvable extension of $K$.  The geometric observation that enables us
to find such points is the following.  Let $L/K$ be a line contained
in the variety $Q_1^{\bar{\rho}} \cap Q_2^{\bar{\rho}}$.  Then $L \cap
Q_4$ is a finite subscheme $\Sigma$ of degree $4$ of $S^{\bar{\rho}}$.
Generically, $\Sigma$ will split into $4$ distinct points over a
degree $4$ (whence solvable!) extension of $K$.  Now $Q_1^{\bar{\rho}}
\cap Q_2^{\bar{\rho}}$ is isomorphic to a quadric hypersurface in
$\P^4$, so its Fano variety is rational.  This means we have plenty of
lines in $Q_1^{\bar{\rho}} \cap Q_2^{\bar{\rho}}$, whence plenty of
points in $S^{\bar{\rho}}$ defined over solvable extensions of $K$.
What remains is to make sure we can find such points which satisfy the
local conditions at $3,5,$ and $\infty$ required in the theorem.  Our
strategy will be to define suitable lines over completions of $K$ at
the relevant primes, and finally to use strong approximation on the
Fano variety $F_1(Q_1^{\bar{\rho}} \cap Q_2^{\bar{\rho}})$  to find a
global line which is adelically close to the specified local ones.

\subsection{Archimedean places}
 
Let $c$ be a complex
conjugation in $\Gal(\bar{K}/K)$, and let $u$ be the corresponding
real place of $K$.

The fact that $\bar{\rho}$ is odd
implies that $\bar{\rho}(c)$ is conjugate to
\beq
\mat{-1}{0}{0}{1}.
\eeq
In particular, if
\beq
\bar{\rho}_0 = \mat{\bar{\chi}_3}{0}{0}{1},
\eeq
we have
\beq
\bar{\rho}_0 | \Gal(\C/K_u) \cong \bar{\rho} | \Gal(\C/K_u),
\eeq
whence
\beq
S^{\bar{\rho}} \times_K K_u \cong 
S^{\bar{\rho}_0} \times_\Q K_u = S^{\bar{\rho}_0} \times_\Q \R.
\eeq

If $s_0, \ldots, s_5$ are our standard coordinates on $S$, we may take
$t_0, \ldots, t_5$ as coordinates on $S^{\bar{\rho}_0}_\Q$, where
\begin{equation}
(s_0, \ldots, s_5) = (t_0 + \sqrt{-3}t_1, t_0 - \sqrt{-3}t_1,
t_2 + \sqrt{-3}t_3, t_2 - \sqrt{-3}t_3,t_4,t_5).
\label{e:barrhocoords}
\end{equation}

Then $S^{\bar{\rho}_0}_\Q$ is isomorphic to the
complete intersection
\begin{eqnarray*}
&&\sigma_1(t_0 + \sqrt{-3}t_1, t_0 - \sqrt{-3}t_1,
t_2 + \sqrt{-3}t_3, t_2 - \sqrt{-3}t_3,t_4,t_5) \\ & = &
\sigma_2(t_0 + \sqrt{-3}t_1, t_0 - \sqrt{-3}t_1,
t_2 + \sqrt{-3}t_3, t_2 - \sqrt{-3}t_3,t_4,t_5) \\
& = & \sigma_4(t_0 + \sqrt{-3}t_1, t_0 - \sqrt{-3}t_1
t_2 + \sqrt{-3}t_3, t_2 - \sqrt{-3}t_3,t_4,t_5) \\
& = & 0.
\end{eqnarray*}

Now choose a real line $L_\R$ in $F_1(Q_1^{\bar{\rho}_0} \cap
Q_2^{\bar{\rho}_0})(\R)$ with the property that $L_\R \cap
S^{\bar{\rho}_0}$ consists of four distinct real points.  For
instance, we may choose $L_\R$ to be the line 
\beq
(t_0,t_1,t_2,t_3,t_4,t_5) = (7/15 + (4/3)t, -1, (4/5)-t, t, -2,
-8/15-(2/3)t).
\eeq

Let $L_u$ be the corresponding line in $F_1(Q_1^{\bar{\rho}} \cap
Q_2^{\bar{\rho}})(K_u)$.

\subsection{Primes above $5$}

\label{ss:char5}

Let $K_v$ be the completion of $K$ at a prime $v$ dividing $5$, and
let $E_0$ be the splitting field of $\bar{\rho}|G_{K_v}$.  Note that,
by hypothesis, $E_{v'}$ has odd absolute ramification degree.

As above, our aim is to find a suitable line in $(Q_1^{\bar{\rho}}
\cap Q_2^{\bar{\rho}})$ over some unramified extension of $E_{v'}$.  Since
$\bar{\rho}$ is trivial on $\Gal(\bar{\Q}_5/E_{v'})$, we may drop
superscripts and look for lines on $Q_1 \cap Q_2$.

\begin{lem} There exists a finite unramified extension $E$ of $E_{v'}$ and
a line $L_v/E$ contained in $Q_1 \cap Q_2/E$ such that
\begin{itemize}
\item $L_v$ is disjoint from the cuspidal locus $Z$;
\item $(L_v \cap Q_4)(E)$ consists of $4$ distinct $E$-points;
\item For each $x \in (L_v \cap Q_4)(E)$, the functions
\beq
\sigma_5^{-6}(\sigma_3^2 - 4\sigma_6)^5
\eeq
and
\beq
\sigma_5^{-3}\sigma_3^{-1}(\sigma_3^2 - 4\sigma_6)^3
\eeq
have non-positive valuation when evaluated at $x$.
\end{itemize}
\end{lem}

\begin{proof}
One checks that $Q_1 \cap Q_2$ is isomorphic over $\Z_5^{unr}$ to the
Pl\"{u}cker quadric threefold $T := V(y_0 y_1 + y_2 y_3 + y_4^2) \subset
\P^4$.  We also know (see \cite[\S 6, Ex.\ 22.6]{harr:ag}) an explicit
$3$-parameter family of lines on $T$, which is to say a map
\beq
\lambda: \P^3/\Spec \Z_5 \ra F_1(T);
\eeq
moreover, $\lambda$ is an isomorphism over any algebraically closed
field.  Composing $\lambda$ with an isomorphism between $T$ and $Q_1
\cap Q_2$ yields a map
\beq
L:  \P^3/\Spec \Z_5^{unr} \ra F_1(Q_1 \cap Q_2)
\eeq
which is an isomorphism over any algebraically closed field.

The set of $\bar{\mathbf{p}} \in \P^3(\fivebar)$ such that $L({\bar{\mathbf{p}}}) \cap
Q_4/\fivebar$ consists of $4$ distinct $\fivebar$-points is
Zariski-open.  To check that it is not empty, we need only exhibit
a single such line $L$ in $(Q_1 \cap Q_2)/\fivebar.$   One such line is
\beq
(s_0,s_1,s_2,s_3,s_4,s_5)  = ((1-\sqrt{-3})t, (1+\sqrt{-3})t, -t +
(1+\sqrt{-3})u, -t + (1-\sqrt{-3})u, t, -t-2u).
\eeq
One checks that the restriction of $Q_4$ to $L$ is $-3t(8u^3 - t^3)$,
which indeed has $4$ distinct roots over $\fivebar$.

Let $V$ be the closed
subscheme of $S/\fivebar$ where the form $\sigma_3^2 - 4\sigma_6$
vanishes.  Then $V$ is a curve.  Moreover, if $x$ is a point in
$S/\fivebar$, the subscheme of $\P^3/\fivebar$ parametrizing lines passing
through $x$ is one-dimensional.  So the subscheme of $\P^3/\fivebar$
parametrizing lines intersecting $V$ is at most two-dimensional.  We
may thus choose a point $\bar{\mathbf{p}} \in \P^3(\fivebar)$ such that
$L({\bar{\mathbf{p}}}) \cap Q_4/\fivebar$ consists of four distinct
$\fivebar$-points, none contained in $V$.

Now let $\mathbf{p}$ be a lift of $\bar{\mathbf{p}}$ to $\P^3(\Q_5^{nr})$.  Then
$L(\mathbf{p})$ is a line contained in $Q_1 \cap Q_2$ whose
intersection with $Q_4$ consists of four distinct points defined over
some unramified extension of $\Q_5$.  Let $E$ be the compositum of
this extension with $E_{v'}$.  Since $Z$ is one-dimensional, we may choose
$\mathbf{p}$ such that $L(\mathbf{p}) \cap Q_4$ is disjoint from $Z$,
by the same argument as above.

Let $x$ be a point in $L(\mathbf{p}) \cap Q_4$, and choose integral coordinates
for $x$ so that at least one coordinate has non-positive valuation.
Then $(\sigma_3^2-4\sigma_6)(x)$ has
non-positive valuation, so the third desired condition on $L(\mathbf{p})$ is
satisfied.  This completes the proof.
\end{proof}

Now take $L$ and $E$ as in the lemma.  Let $x_1, x_2, x_3, x_4$
be the four $E$-points making up $(L_v \cap S)(E)$.  Then
each $x_i$ corresponds to an abelian variety $A_i/E$ with real
multiplication by $\OO$ admitting an isomorphism $A[3]
\cong \bar{\rho} \cong \F_9^{\oplus 2}$ of $\OO$-module schemes over $E$.
Since $\XX^{\bar{\rho}}/\OO_E$ is a fine proper moduli space, we know that
$A_i$ extends to a semi-abelian scheme $\AA_i/\OO_E$.

We now want to show that each $A_i$ has good ordinary or
multiplicative reduction.  We begin by recalling some facts on Hilbert modular
forms for $\SL_2(\OO)$.  By a result of Nagaoka~\cite[Th.\ 2]{naga:hmfz}, the ring
$M_{2*}(\SL_2(\OO),\Z[1/2])$ of even-weight symmetric modular forms over $\Z[1/2]$ for
this group is generated by forms
$\phi_2, \chi_6,$ and  $\chi_{10}$ of weights
$2,6,10$. The form $\phi_2$ is the weight $2$
Eisenstein series, while $\chi_6$ and $\chi_{10}$ are cuspforms.

By restricting to various modular curves
on $X$, and comparing $q$-expansions, one can verify the following
identities between the level $1$ forms above and the symmetric
functions $\sigma_i$ in $s_0, \ldots, s_5$:
\begin{eqnarray*}
\phi_2 & = & -3 \sigma_5^{-1} (\sigma_3^2 - 4 \sigma_6) \\
\chi_6 & = & \sigma_3 \\
\chi_{10} & = & (-1/3)\sigma_5.
\end{eqnarray*}

(Note that the constants here depend on our original choice of the
weight $2$ forms $s_i$.  Modifying that choice by a constant $c$ would
modify each formula above by $c^{k/2}$, where $k$ is the weight of the
modular form in the expression.)

So by our choice of $L_v$, the modular functions $\phi_2^3/\chi_6$ and
$\phi_2^5/\chi_{10}$ have non-positive valuation when evaluated on
$A_i$.  The desired ordinarity now follows from the following lemma.

\begin{lem} Let $A/E$ be an semi-HBAV over a finite extension
$\OO_E/\Z_5$.  Suppose that the modular functions $\phi_2^3/\chi_6$ and
$\phi_2^5/\chi_{10}$ evaluated at $A$ have non-positive valuation.
Then $A$ has potentially ordinary or multiplicative reduction.
\label{le:ord5}
\end{lem}

\begin{proof}

Let $\Omega$ be the determinant of the pushforward of the relative
cotangent sheaf of $A/\OO_E$.  Then $\Omega$ is a free rank $1$
$\OO_E$-module.  Let $\omega$ be a section generating $\Omega$.  Then
every modular form $f$ with coefficients $\OO_E$ has a well-defined
value $f(A,\omega)$.  Suppose $\phi_2(A,\omega) \in \ic{m}_E$.  Then
by the hypothesis of the theorem, we have also that $\chi_6(A,\omega)$
and $\chi_{10}(A,\omega) \in \ic{m}_E$.  But this is impossible, as we
show in the following paragraph.

Let $q \geq 7$ be a prime, let $\iota:(\Z/q\Z)^2 \cong A[q]$ be an arbitrary
full level $q$ structure on $A$, and let $f$ be a modular form of
weight $2k$ and full level $q$.  Since every even-weight modular form
of full level $q$ is integral over $M_{2*}(\Gamma(1),\Z) =
\Z_5[\phi_2,\chi_6,\chi_{10}],$ we have $f(A,\iota,\omega) \in
\ic{m}_E$.  But this is impossible, because for $2k$ sufficiently
large, the sheaf $\Omega^{\tensor 2k}$ is very ample on the level $q$ moduli
scheme $\XX(q)$.

We conclude that $\phi_2(A,\omega) \neq \ic{m}_E$.  So the mod $5$
reduction $\bar{\phi}_2(\bar{A},\bar{\omega})$ is not equal to $0$.
Since the reduction mod $5$ of $\phi_2^2$ is the Hasse invariant, $A$
has good ordinary or multiplicative reduction by~\cite{andr:ango}.
\end{proof}

In order to keep straight the primes of $K$ over $5$, we now give the extension
$E/E_{v'}$ the new name $E'_{v'}/E_{v'}$.

\subsection{Primes above $3$}

\label{ss:char3}

We begin by observing that we can apply to $\bar{\rho}$ a global
quadratic twist $\psi$ such that
\begin{equation}
\bar{\rho}|D_w \cong \mat{\bar{\chi}_3}{*}{0}{1}
\label{eq:goodtwist}
\end{equation}
for all primes $w$ of $K$ dividing $3$.  Since twisting does not
affect modularity, we may assume that \eqref{eq:goodtwist} holds for
our original $\bar{\rho}$.

Let $w$ be a prime of $K$ dividing $3$, and let $K_w$ be the
completion of $K$ and $w$.  Now the $*$ in \eqref{eq:goodtwist} is
a cocycle corresponding to an element $\lambda \in K_w^*
\tensor_\Z \F_9$.  Write $\bar{\rho}_\lambda$ for the representation
$\bar{\rho}|G_w$.

\begin{lem}
There exists a line $L_w$ in $\P^5_{K_w}$ satisfying the following
conditions.
\begin{itemize}
\item $L_w$ is contained in $Q_1^{\bar{\rho}_\lambda} \cap
Q_2^{\bar{\rho}_\lambda}$.
\item The intersection $L_w \cap Q_4^{\bar{\rho}_\lambda}$ splits into four
distinct points over an unramified extension $E'_{w'}$ of $K_w$.
\item The four HBAV's $A_1,A_2,A_3,A_4$ corresponding to the four
points of $L_w \cap Q_4^{\bar{\rho}_\lambda}(\Q^{nr}_3)$ have potentially
multiplicative reduction. 
\end{itemize}
\end{lem}

\begin{proof}
We define new coordinates $x_0, \ldots, x_5$ on $\P^5$ by the rule
\begin{eqnarray*}
x_0 & = & \omega s_0 + \omega^2 s_1 + s_4 \\
x_1 & = & \omega^2 s_0 + \omega s_1 + s_4 \\
x_4 & = & s_0 + s_1 + s_4 \\
x_2 & = & \omega s_2 + \omega^2 s_3 + s_5 \\
x_3 & = & \omega^2 s_2 + \omega s_3 + s_5 \\
x_5 & = & s_2 + s_3 + s_5 \\
\end{eqnarray*}
where $\omega$ is a cube root of unity.
With these coordinates, one checks that $Q_1^{\bar{\rho}_\lambda}$ is
defined by $x_4 + x_5$ and $Q_2^{\bar{\rho}_\lambda}$ by
\beq
x_4^2 + x_5^2 - x_0x_1 - x_2x_3 + 3x_4x_5.
\eeq
So a family of lines in $Q_1^{\bar{\rho}_\lambda} \cap
Q_2^{\bar{\rho}_\lambda}$ is given by
\beq
L_{a,b,c}:  x_0 = ax_2 + bx_4, x_3 = -ax_1 + cx_4, x_4 = -(bx_1 + cx_2), x_5 =
-x_4.
\eeq
One checks that the equation for $Q_4^{\bar{\rho}_\lambda}$ is given
by
\begin{eqnarray*}
-3x_0x_1x4x_5 - 3 x_2x_3x_4x_5 + 3x_0x_1x_2x_3 + x_4x_5(x_4^2 +
3x_4x_5 + x_5^2) - 3x_0x_1x_5^2 - 3x_2x_3x_4^2 \\ +
\lambda_1^3 x_0^3x_5 + \lambda_1^{-3} x_1^3x_5
+ \lambda_2^3 x_2^3 x_4 + \lambda_2^{-3} x_3^3 x_4 & &
\end{eqnarray*}
where $\lambda_1,\lambda_2$ are elements of $K$ whose class in
$K^*/(K^*)^3$ is determined by $\lambda$.

The equation for $Q_4^{\bar{\rho}_\lambda}$ restricted to $L_{a,b,c}$
is of the form 
\beq
P = \sum_{i=0}^4 P_i(a,b,c)x_1^i x_2^{4-i}.
\eeq
Suppose that $\ord_w(b)$ and $\ord_w(c)$ are approximately equal and that
both are much greater than $\ord_w(a)$, which is in turn much greater
than $0$.  Then one checks that
\begin{eqnarray*}
P_4(a,b,c) & = & \lambda_1^{-3}b + \mbox{higher order terms} \\
P_3(a,b,c) & = & \lambda_1^{-3}c + \mbox{higher order terms} \\
P_2(a,b,c) & = & -3a^2 + \mbox{higher order terms} \\
P_1(a,b,c) & = & -\lambda_2^3 b + \mbox{higher order terms} \\
P_0(a,b,c) & = & -\lambda_2^3 c + \mbox{higher order terms}.
\end{eqnarray*}

It is an easy calculation that $P$ then factors over $K^{unr}_w$ into
a constant and two quadratics, one of which has discriminant with
valuation equal to that of $\lambda_1^{-3} a^{-2} b$, the other of which has
discriminant with valuation equal to that of $\lambda_2^3 a^{-2} c$.
In particular, if $b$ and $c$ are chosen so that $\ord_w(b)$ has the
same parity as $\ord_w(\lambda_1)$, and $\ord_w(c)$ the same parity as
$\ord_w(\lambda_2)$, the points of $L_{a,b,c} \cap
Q_4^{\bar{\rho}_\lambda}$ are defined over an unramified extension of $K_w$.

Note that since $\ord_w(a), \ord_w(b),$ and $\ord_w(c)$ are much
larger than $0$, we know that $L_{a,b,c}$ is very $v$-adically close
to the cusp line $L_{0,0,0}$.  Moreover, the form for $P$ given above
shows that the four points in $L_{a,b,c} \cap
Q_4^{\bar{\rho}_\lambda}$ are $w$-adically close to $[0:1:0:0:0:0]$
and $[0:0:1:0:0:0]$.  These points do not lie in the image of any
non-cuspidal point of the Hilbert modular surface
$X$~\cite[VIII(2.6)]{vdge:hms}.  So the abelian varieties $A_1, A_2,
A_3, A_4$ corresponding to these four points  to the points of $L_{a,b,c} \cap
Q_4^{\bar{\rho}_\lambda}$ have potentially multiplicative reduction at
$w$, as desired. 

\end{proof}

The fact that $A_i$ has potentially multiplicative reduction implies
that the $\ell$-adic representation $\rho_{A_i,\ell}$ attached to
$A_i$ satisfies  
\beq
\rho_{A_i,\ell}|I_{w'} \cong \mat{*}{*}{0}{*}.
\eeq
The diagonal entries are finite-order characters with image in
$\OO^*$, which must have image in $\pm 1$.  It follows that either
$A_i$ or its quadratic twist $A_i \tensor \bar{\chi}_3$ has
semi-stable reduction at $w'$. Note that
$A_i[3](\bar{\Q}_3)$ contains a canonical subgroup $A^0_i$, consisting
of those points reducing to the identity in the N\'{e}ron model of
$A_i$.  Now the action of $I_{w'}$ on $A_i[3]$ is isomorphic to
\beq
\mat{\bar{\chi}_3}{*}{0}{1},
\eeq
and $A_i[3]$ has semistable reduction if and only if the canonical
subgroup is the subgroup on which $I_{w'}$ acts as $\bar{\chi}_3$.
Note that $A_i[3]$ is semistable whenever $*$ is nonzero.

Each of the points in $L_{a,b,c} \cap Q_4^{\bar{\rho}_\lambda}$
corresponds to an abelian variety $A_i$ together with an isomorphism
$\phi_i$ from $A_i[3]$ to $V_{\bar{\rho}_\lambda}$.  The
$(A_i,\phi_i)$ correspond to points of $X(\Q_3)$ lying in a small
$3$-adic neighborhood of a cusp; since the canonical subgroup varies
$3$-adically continuously, the image $\phi_i(A^0_i)$ is the same subgroup of
$V_{\bar{\rho}_\lambda}$ for each $i$. In particular, either all four
of the $A_i$ have semistable reduction at $w'$, or all four have
semistable reduction after twisting by $\bar{\chi}_3$.

\subsection{The global construction}

We now combine the local arguments above into the global statement we
desire.

Choose a finite Galois extension $K'/K$ such that;
\begin{itemize}
\item $K'$ is totally real;
\item $K'/K$ is solvable;
\item The completion of $K'$ at any prime $v$ above $5$ is isomorphic to
an unramified extension of $E'_{v'}$;
\item The completion of $K'$ at any prime $w$ above $3$ is isomorphic to
an unramified extension of $E'_{w'}$;
\item  $F_1((Q_1 \cap Q_2)^{\bar{\rho}})$ is rational over $K'$.
(Since $F_1((Q_1 \cap Q_2)^{\bar{\rho}})$ is geometrically rational,
this amounts to trivializing an element of the Brauer group.)
\end{itemize}

(See \cite[Lemma 2.2]{tayl:icos2} for the existence of $K'$.)

Since $F_1((Q_1 \cap Q_2)^{\bar{\rho}})$ is a rational variety, we can choose $L \in
F_1((Q_1 \cap Q_2)^{\bar{\rho}})(K')$ such that the image of $L$ under the map
\beq
F_1((Q_1 \cap Q_2)^{\bar{\rho}})(K') \ra \bigoplus_{v_i | 5}
F_1((Q_1 \cap Q_2)^{\bar{\rho}})(K'_{v_i})
\oplus 
\bigoplus_{w_i | 3}
F_1((Q_1 \cap Q_2)^{\bar{\rho}})(K'_{w_i})
\bigoplus_{u | \infty} \F_1((Q_1 \cap Q_2)^{\bar{\rho}})(K'_{u})
\eeq
is arbitrarily adelically close to $(L_{v_1}, \ldots, L_{w_1}, \dots,
L_{u_1},\ldots)$.

The intersection $L \cap S^{\bar{\rho}}$ is a zero-dimensional scheme
of degree $4$ over $K'$.  Modifying our choice of $L$ if necessary, we
can arrange for $L \cap S$ to be in the image of the rational map from
$X^{\bar{\rho}}$.  Let $F$ be a splitting field for $L \cap
S^{\bar{\rho}}$.  Note that $F$ is solvable over $K'$, whence also
over $K$.  Then we can think of $L \cap S$ as specifying four HBAV's
$A_i/F$, with $A_i[3]/F \cong V_{\bar{\rho}}/F$.
 
By our choices of $L_u$, the field $F$ is totally
real.  Similarly, our choices of $L_v$ and $L_w$ guarantee that $A_i$
and $F$ satisfy the local conditions at $3$ and $5$ stated in the
theorem, by a theorem of Kisin~\cite{kisi:zeit}.

It remains only to check that $L$ can be chosen so that
$A_i[\sqrt{5}]$ is an absolutely irreducible $\GalF$-module, for some
$i$.  Let $\pi: T_\C \ra F_1(Q_1 \cap Q_2)_\C$ be the degree
$4$ cover defined by $\pi^{-1}(L) = (L \cap S)(\C)$.  Then $T$
projects to $S$, and the covering of $S$ by the Hilbert modular
surface for the congruence subgroup of full level $\sqrt{5}$ pulls
back to a covering $T' \ra T$.  The covering $T'/T$ is Galois with
group $\PSL_2(\F_5)$, and the Galois group of the cover $T'/F_1(Q_1
\cap Q_2)$ is thus a subgroup $G$ of the wreath product of $S_4$ and
$\PSL_2(\F_5)$.

We want to show that $G$ is large.  Let $B$ be an HBAV over a number
field $M$ such that the
map
\beq
\Gal(\bar{\Q}/M) \ra \GL(B[3]) \oplus \GL(B[\sqrt{5}])
\eeq
is surjective, and let $p$ be a point of $S(M(B[3]))$ corresponding to
$B$.  Let $\ell$ be a line in $Q_1 \cap Q_2$ passing through $p$; we can
choose $\ell$ to be defined over a quadratic extension of $M$.  Now the
fact that the map
\beq
\Gal(\bar{\Q}/M(B[3])) \ra \GL(B[\sqrt{5}])
\eeq
is surjective implies that the Galois group $G$ contains a copy of
$\PSL_2(\F_5)$.  By Ekedahl's version of Hilbert's irreducibility
theorem~\cite{eked:hit}, we can now choose $L$ in the adelic
neighborhood specified above in such a way that one of the four
abelian varieties $A_i/F$ has $A_i[\sqrt{5}]$ an absolutely
irreducible $\GalF$-module.  This completes the proof of
Proposition~\ref{pr:hbav}. 

\section{Modularity}

Now that we have exhibited $\bar{\rho}$ as a representation appearing
on the torsion points of an abelian variety, we can prove that
$\bar{\rho}$ is modular.  Our argument proceeds exactly along the lines of
\cite{tayl:icos2} and \cite{khar:hmfp}.  We begin by recording the
case we need of a theorem of Skinner and Wiles.

\begin{thm}
Let $K$ be a totally real number field, let $p>2$ be a rational prime,
let $L$ be a finite extension
of $\Qp$, and let
\beq
\rho: \GalK \ra \GL_2(L)
\eeq
be a continuous irreducible representation ramified at only finitely many primes.  Suppose
\begin{itemize}
\item $\det \rho = \psi \chi_p^{k-1}$ for some finite-order character
$\psi$;
\item $\rho$ is {\em ordinary} in the sense that, for each prime
$v$ of $K$ dividing $p$,
\beq
\rho|D_v \cong \mat{\chi_p}{*}{0}{1};
\eeq

\item $\rho$ is $D_v$-distinguished for all primes $v$ of $K$ dividing
$p$;
\item There exists an ordinary modular Galois representation $\rho'$
and an isomorphism between the mod $p$ representations $\bar{\rho}$
and $\bar{\rho}'$.
\end{itemize}

Then $\rho$ is modular.
\label{th:sw2}
\end{thm}

\begin{proof}See \cite{skin:sw2}. \end{proof}

We are now ready to prove our main result.

\begin{thm}
Let
\beq
\bar{\rho}: \GalQ \ra \GL_2(\F_9)
\eeq
be an odd, absolutely irreducible Galois representation such that
\begin{itemize}
\item  The restriction of $\bar{\rho}$ to $D_3$ can be written as
\beq
\bar{\rho}|D_3 \cong \mat{\psi_1}{*}{0}{\psi_2},
\eeq
where $\psi_1$ and $\psi_2$ are distinct characters such that $(\psi_1
\oplus \psi_2)|I_3$ is isomorphic to the direct sum of the trivial
character and $\bar{\chi}_3$.

\item The image of the inertia group $I_5$ lies in $\SL_2(\F_9)$, and
has odd order.

\end{itemize}

Then $\bar{\rho}$ is modular.
\label{th:main}
\end{thm}

\begin{proof}

First of all, $\theta = (\det \bar{\rho})^{-1} \bar{\chi}_3$ is a character of
Galois which annihilates complex conjugation, since $\bar{\rho}$ is
odd.  We thus have a totally real abelian extension $K/\Q$ defined by
$\GalK = \ker \theta$.  Since $\det(\bar{\rho})(I_5)$ is trivial,
$I_5$ lies in the kernel of $\theta$, and $K$ is unramified at $5$.
Likewise, $\theta(I_3)$ is trivial, so $K$ is unramified at $3$.

By applying a quadratic twist, we may assume that
\beq
\bar{\rho}|D_3 \cong \mat{\bar{\chi}_3}{*}{0}{1}.
\eeq

Now the conditions on $\bar{\rho}|D_3$ and $\bar{\rho}|D_5$ imply the
corresponding local conditions in Proposition~\ref{pr:hbav}.  We may
now choose an extension $F_0/K$ and an abelian variety $A/F_0$ satisfying
the four hypotheses given in that Proposition.  

From here, we proceed along the lines of \cite{tayl:icos2}.  First,
we claim that the irreducible representation
\beq
\bar{\rho}_{A,\sqrt{5}}: \Gal(\bar{\Q}/F_0) \ra \GL_2(\F_5)
\eeq
induced by the torsion subscheme $A[\sqrt{5}]$ is modular.

Now it follows from the discussion in section~\ref{ss:char3} that
either $A$ or $A \tensor \bar{\chi}_3$ has multiplicative reduction at all
primes of $F_0$ over $3$.  In case it is $A \tensor \bar{\chi}_3$
which is semistable, return to the beginning, replace $\bar{\rho}$ by
$\bar{\rho} \tensor \bar{\chi}_3$, and start over.  We may now assume
$A$ has multiplicative reduction at all primes of $F_0$ over $3$.

Now the subgroup
\beq
\bar{\rho}_{A,\sqrt{5}}(I_w) \subset \SL_2(\F_5)
\eeq
is unipotent.  Thus, we can find a totally real solvable extension
$F/F_0$, unramified over $5$ and with odd ramification
degree at every prime over $3$,  such that $D_w$ acts trivially on
$A[\sqrt{5}]$ for every prime $w$ of $F$ dividing $3$.  Then the twist
of the modular curve $X(5)_{F}$ by $\bar{\rho}_{A,\sqrt{5}}$ is
isomorphic to $X(5)$ when base changed to any $3$-adic completion of
$F$.  In particular, there exists an elliptic curve $E/F$ such that $E[5]/F \cong A[\sqrt{5}]/F$ and $E$
has good ordinary reduction at each prime of $F$ over $3$.  By
another use of Ekedahl~\cite{eked:hit}, we can assume that $E[3]$ is
an absolutely irreducible $\Gal(\bar{\Q}/F)$-module.

Since $F$ has odd absolute ramification degree at all primes over
$3$, we have by Proposition~\ref{pr:distinguished} that $E[3]/F$ is
distinguished at all primes over $3$.  The mod-$3$ representation
$E[3]$ is modular by the Langlands-Tunnell theorem, and the $3$-adic
Galois representation $T_3 E$ is modular by Theorem~\ref{th:sw2}.  It
follows that $T_5 E$ is modular, and so $E[5]$, whence also
$A[\sqrt{5}]/F$, is modular.

By hypothesis, $A$ has good ordinary or multiplicative
reduction at $5$, so $T_{\sqrt{5}} A$ is an ordinary
representation.  Because $F/\Q$ has odd ramification degree over $5$,
$T_{\sqrt{5}} A$ is $D_v$-distinguished for all primes $v$ dividing
$5$.  Now $T_{\sqrt{5}} A / F$ is modular by another application of
Theorem~\ref{th:sw2}.  This implies that $T_3 A/F$ is also modular,
which in turn implies the modularity of $A[3]/F$, which
is the restriction to $\GalF$ of our original representation $\bar{\rho}$.

Recall that  the restriction of $\bar{\rho}$ to the decomposition group
$D_3$ is of the form
\beq
\bar{\rho}|D_3 \cong \mat{\bar{\chi_3}}{*}{0}{1}.
\eeq

 By results of Ramakrishna as refined by Taylor~\cite[Thm
 1.3]{tayl:icos2}, $\bar{\rho}$ can be lifted to a $3$-adic representation   
 
\beq
\rho: \GalQ \ra \GL_2(W(\F_9))
\eeq
such that
\beq
\rho|D_3 \cong \mat{\chi_3} {*}{0}{1}.
\eeq

Now $\rho|\GalF$ is an ordinary, $D_3$-distinguished $3$-adic
representation of $\GalF$ whose reduction mod $3$ is isomorphic to the
modular representation $A[3]$.  Applying Theorem~\ref{th:sw2}
once more, using $T_3 A$ as $\rho'$, we have that $\rho|\GalF$ is modular.  Now we argue by
cyclic descent as in \cite{tayl:icos2}.  Let $F'$ be a subfield of $F$ such that $F/F'$
is a cyclic Galois extension.  Then the automorphic form $\pi$ on $\GL_2(F)$
corresponding to $\rho$ is preserved by $\Gal(F/F')$.  Therefore,
$\pi$ descends to an automorphic form on $\GL_2(F')$.  Continuing
inductively, one finds that $\rho$ itself is associated to a modular
form on $\GL_2(\Q)$; therefore, its mod $3$
reduction $\bar{\rho}$ is modular. 
\end{proof}

This case of Serre's conjecture can be used to prove the modularity of
Hilbert-Blumenthal abelian surfaces, under some conditions on
reduction at $3$ and $5$. 

\begin{cor}  Let $A/\Q$ be a Hilbert-Blumenthal abelian surface which has
good ordinary reduction at $3$ and $5$.  Then $A$ is a quotient of
$J_0(N)$ for some integer $N$. 
\end{cor}

\begin{proof}  Let $v$ be a prime of the field of real multiplication
dividing $3$.  If $A[v]$ is absolutely reducible, then the corollary
follows from a theorem of Skinner and Wiles~\cite{skin:sw1}.  So we
may assume that $A[v]$ is absolutely irreducible.

If $3$ is split or ramified in the ring of real multiplication $\OO$,
then the corollary follows from Langlands-Tunnell applied to $A[v]$
followed by Diamond's refinement of the theorem of Wiles and
Taylor-Wiles~\cite{diam:hecke}.  If, on the other hand, $3$ is inert in $A$, then $A[3]$
yields a representation $\rho_{A,3}: \GalQ \ra \GL_2(\F_9)$, which is
easily seen to satisfy the conditions of Theorem~\ref{th:main}.  So
$A[3]$ is modular, and it follows, again by Diamond's theorem, that
$A$ is modular.
\end{proof}

\end{document}